\documentclass[11pt]{amsart}
\usepackage[T1]{fontenc}
\usepackage{times}

\makeatletter

\newcommand{\LyX}{L\kern-.1667em\lower.25em\hbox{Y}\kern-.125emX\spacefactor1000}

\theoremstyle{plain}    
\newtheorem{thm}{Theorem}[section]
\numberwithin{equation}{section} %% Comment out for sequentially-numbered
\numberwithin{figure}{section} %% Comment out for sequentially-numbered
\theoremstyle{plain}    
\newtheorem{cor}[thm]{Corollary} %%Delete [thm] to re-start numbering
\theoremstyle{plain}    
\newtheorem{lem}[thm]{Lemma} %%Delete [thm] to re-start numbering
\theoremstyle{plain}    
\newtheorem{prop}[thm]{Proposition} %%Delete [thm] to re-start numbering
\theoremstyle{remark}    
\newtheorem*{acknowledgement*}{Acknowledgement} 

\newcommand{\sa}{\operatorname{sa}}
\usepackage{amssymb}
\makeatother

\begin{document}

\title{Free Fisher Information for Non-Tracial States}

\author{Dimitri Shlyakhtenko}

\date{\today}

\address{Department of Mathematics, UCLA, Los Angeles, CA 90095}

\email{shlyakht@math.ucla.edu}

\begin{abstract}
We extend Voiculescu's microstates-free definitions of free Fisher information
and free entropy to the non-tracial framework. We explain the connection between
these quantities and free entropy with respect to certain completely positive
maps acting on the core of the non-tracial non-commutative probability space.
We give a condition on free Fisher information of an infinite family of variables,
which guarantees factoriality of the von Neumann algebra they generate.
\end{abstract}
\maketitle

\newcommand{\reals}{\mathbb {R}}

\newcommand{\bbc}{\mathbb {C}}

\section{Introduction.}

Free entropy and free Fisher information were introduced by Voiculescu \cite{dvv:entropy1},
\cite{dvv:entropy2}, \cite{dvv:entropy5} in the context of his free probability
theory \cite{DVV:book} as analogs of the corresponding classical quantities.
These quantities are usually considered in the framework of tracial non-commutative
probability spaces; not surprisingly, the most striking applications of free
entropy theory were to tracial von Neumann algebras (see e.g. \cite{dvv:entropy3},
\cite{ge:entropy2}, \cite{stephan:thinness}). Recently, however, it turned
out that some type III factors associated with free probability theory \cite{shlyakht:quasifree:big}
have certain properties in common with their type II\( _{1} \) cousins \cite{shlyakht:prime}.
This gives rise to a speculation that there is room for free entropy to exist
outside of the context of tracial non-commutative probability spaces.

The goal of this paper is to initiate the development of free Fisher information,
based on Voiculescu's microstates-free approach \cite{dvv:entropy5}, in the
non-tracial framework. The key idea is that all of the ingredients going into
the definition of free Fisher information in this case must behave covariantly
with respect to the modular group \cite{stratila:modular} of the non-tracial
state. The principal example of a family of variables for which free Fisher
information is non-trivial, and which belong to an algebra not having any traces,
are semicircular generators of free Araki-Woods factors, taken with free quasi-free
states \cite{shlyakht:quasifree:big}.

We describe another route towards free Fisher information, which is based on
first converting the non-tracial von Neumann algebra into a larger algebra,
the core (having an infinite trace), and then considering free Fisher information
relative to a certain completely positive map (in the spirit of \cite{shlyakht:cpentropy}).
We should point out that it is this approach that is most likely to connect
with the microstates free entropy (as suggested by \cite{guionnet:band}; see
also \cite{shlyakht:relmicro}), since it is at present unclear what a microstates
approach to free entropy in the non-tracial framework should be.

We finish the paper with a look at free Fisher information on von Neumann algebras
that have traces. Our first result is that once the algebra has a trace, the
free Fisher information is automatically infinite when computed with respect
to a non-tracial state. It is likely that on any von Neumann algebra, free Fisher
information can be finite for only very special states (however, we do not have
any results in this direction in the non-tracial category). Another result of
the present paper is a statement guaranteeing factoriality of a tracial von
Neumann algebra, once we know that it has an infinite generating family whose
free Fisher information is bounded in a certain way.

\begin{acknowledgement*}
This work was supported by an NSF postdoctoral fellowship. The author is grateful
to MSRI and the organizers of the operator algebras program for their hospitality.
\end{acknowledgement*}

\section{Free Fisher Information for Arbitrary KMS States.}

\subsection{Free Brownian motion in the presence of a modular group.}

Let \( M \) be a von Neumann algebra, \( \phi :M\to \mathbb {C} \) be a normal
faithful state on \( M \). Denote by \( \sigma ^{\phi }_{t} \) the modular
group of \( \phi  \). Denote by \( L^{2}(M^{\sa },\phi )\subset L^{2}(M,\phi ) \)
the closure of the real subspace of self-adjoint elements \( M^{\sa }\subset M \). 

Let \( X=X^{*}\in M \) and let \( B\subset M \) be a subalgebra. Assume that
\( \sigma _{t}^{\phi }(B)\subset B \) for all \( t\in \reals  \). 

Consider the von Neumann algebra \( \mathcal{M}=\Gamma (L^{2}(M^{\sa },\phi )\subset L^{2}(M,\phi )) \),
taken with the free quasi-free state \( \phi _{\mathcal{M}} \). Consider the
element \( Y=s(X)\in \mathcal{M} \) (see \cite{shlyakht:quasifree:big} for
definitions and notation). Then \( \phi _{\mathcal{M}}(Y\sigma _{t}^{\phi _{\mathcal{M}}}(Y))=\phi (X\sigma _{t}^{\phi }(X)) \),
for all \( t\in \reals  \). 

Consider the algebra \( \mathcal{N}=(M,\phi )*(\mathcal{M},\phi _{\mathcal{M}}) \),
and denote by \( \hat{\phi } \) the free product state on \( \mathcal{N} \).
Note that \( \sigma _{t}^{\hat{\phi }}=\sigma _{t}^{\phi }*\sigma _{t}^{\phi _{\mathcal{M}}} \).
The elements \( X_{\epsilon }=X+\sqrt{\epsilon }Y \), \( \epsilon \geq 0 \)
form a natural free Brownian motion, which behaves nicely under the action of
the modular group. In particular, note that for all \( \epsilon \geq 0 \) and
\( t\in \reals  \),
\[
\hat{\phi }(X_{\epsilon }\sigma ^{\hat{\phi }}_{t}(X_{\epsilon }))=\phi (X\sigma _{t}^{\phi }(X))\cdot (1+\epsilon ).\]
Furthermore, for each \( \epsilon >0 \), the distribution of \( X_{\epsilon } \)
is that of a free Brownian motion at time \( \epsilon  \), starting at \( X \);
this is because \( Y \) is a semicircular variable, free from \( X \).

\subsection{Conjugate variables.}

Let \( B[X] \) denote that algebra generated by \( B \) and all translates
\( \sigma _{t}^{\phi }(X) \), \( t\in \reals  \). Assume that \( \{\sigma _{t}^{\phi }(X)\} \)
are algebraically free over \( B \), i.e., satisfy no algebraic relations modulo
\( B \). Denote by \( \partial _{X}:B[X]\to \mathcal{N} \) the derivation
given by:

\begin{enumerate}
\item \( \partial _{X}(\sigma _{t}^{\phi }(X))=\sigma _{t}^{\hat{\phi }}(Y) \)
\item \( \partial _{X}(b)=0 \), \( b\in B \).
\end{enumerate}
Notice that the range of \( \partial _{X} \) actually lies in the subspace
\( B[X]YB[X]\subset \mathcal{N} \). Note also that since \( \partial _{X}(\sigma _{t}^{\phi }(X)) \)
is self-adjoint, we have that for \( P\in B[X] \), \( \partial _{X}(P^{*})=\partial _{X}(P)^{*} \),
i.e., \( \partial _{X} \) is a \( * \)-derivation. Observe finally that \( \partial _{X} \)
is covariant with respect to the modular groups \( \sigma ^{\phi }_{t} \) and
\( \sigma _{t}^{\hat{\phi }} \):
\[
\partial _{X}(\sigma _{t}^{\phi }(P))=\sigma _{t}^{\hat{\phi }}(\partial _{X}(P)),\quad P\in B[X].\]

Define the conjugate variable \( J_{\phi }(X:B)\in L^{2}(B[X],\phi ) \) to
be such a vector \( \xi  \) that
\begin{equation}
\label{eqn:defofJ}
\langle \xi ,P\rangle _{L^{2}(B[X],\phi )}=\langle Y,\partial _{X}(P)\rangle _{L^{2}(\mathcal{N},\hat{\phi })},\quad \forall P\in B[X],
\end{equation}
if a vector \( \xi  \) satisfying such properties exists. Formally, this means
that \( \xi =\partial _{X}^{*}(Y) \), where \( \partial _{X}:L^{2}(B[X],\phi )\to L^{2}(\mathcal{N},\hat{\phi }) \)
is viewed as a densely defined operator.

It is clear, because of the density of \( B[X] \) in \( L^{2}(B[X],\phi ) \),
that \( \xi  \) is unique, if it exists.

It is convenient to talk about \( J_{\phi }(X:B) \) even in the case that \( \{\sigma _{t}(X)\}_{t\in \reals } \)
are not algebraically free over \( B \) (such is the case, for example, when
\( \phi  \) is a trace, and hence \( \sigma _{t}^{\phi }(X)=X \) for all \( t \)).
In this case, one can view \( \partial _{X} \) as a multi-valued map, the set
of values given by the results of application of the definition of \( \partial _{X} \)
in all possible ways; the definition of \( J_{\phi } \) is then that (\ref{eqn:defofJ})
is valid for all values of \( \partial _{X} \).

Note that \( J_{\phi }(X:B) \) depends on more than just the joint distribution
of \( X \) and \( B \) with respect to the state \( \phi  \); it depends
on the joint distribution of the family \( B\cup \{\sigma _{t}^{\phi }(X):t\in \reals ) \).

We continue to denote by \( \sigma _{t}^{\phi } \) the extension of \( \sigma _{t}^{\phi } \)
to the Hilbert space \( L^{2}(M,\phi ) \) (this is precisely the one-parameter
group of unitaries \( \Delta _{\phi }^{it} \), where \( \Delta _{\phi } \)
is the modular operator). In particular, if \( \phi  \) is a trace, then the
definition of \( J_{\phi } \) is (up to a multiple) precisely that of the conjugate
variable of Voiculescu \cite{dvv:entropy5}.

\begin{lem}
Assume that \( \xi =J_{\phi }(X:B) \) exists. Then \( \xi \in L^{2}(M^{\sa },\phi ) \),
and \( \sigma _{t}^{\phi }(J_{\phi }(X:B))=J_{\phi }(\sigma _{t}(X):B) \). 
\end{lem}
\begin{proof}
We note that, because \( \partial _{X} \) is a \( * \)-derivation,
\[
\langle P^{*},\xi \rangle =\overline{\langle Y,\partial _{X}(P^{*})\rangle }=\langle \partial _{X}(P),Y\rangle =\overline{\langle \xi ,P\rangle }.\]
From this it follows that \( \xi  \) is in the domain of the \( S \) operator
of Tomita theory, and moreover that \( S\xi =\xi  \). Hence \( \xi \in L^{2}(M^{\sa }) \). 

One also has
\begin{eqnarray*}
\langle \sigma _{t}^{\phi }(\xi ),P\rangle  & = & \langle Y,\partial _{X}(P)\rangle \\
 & = & \langle \sigma ^{\phi }_{t}(Y),\partial _{\sigma ^{\phi }_{t}(X)}(P)\rangle ,
\end{eqnarray*}
since the joint distributions of \( B[X] \) and \( \{\sigma ^{\hat{\phi }}_{t}(Y)\}_{t\in \reals } \)
is the same as \( B[X] \) and \( \{\sigma _{s+t}^{\hat{\phi }}(Y)\}_{t\in \mathbb {R}} \),
for any \( s \). It follows that \( \sigma _{t}^{\phi }(J_{\phi }(X:B))=J_{\phi }(\sigma _{t}^{\phi }(X):B) \).
\end{proof}
\begin{lem}
\label{lemma:xileftright}Let \( P,Q\in B[X] \), and assume that \( \xi =J_{\phi }(X) \)
exists and is in \( M \). Then
\[
\phi (P\xi Q)=\hat{\phi }(PY\partial _{X}(Q))+\hat{\phi }(\partial _{X}(P)YQ).\]
 
\end{lem}
\begin{proof}
Recall that \( \phi  \) (and \( \hat{\phi } \)) satisfy the KMS condition:
for all \( a,b\in M \) (or \( \in \mathcal{N} \)), there exists a (unique)
function \( f(z) \), analytic on the strip \( \{z:0<\Im z<1\} \), and so that
(writing \( \sigma _{t} \) for either \( \sigma ^{\phi }_{t} \) or \( \sigma ^{\hat{\phi }}_{t} \)
)
\begin{eqnarray*}
\phi (a\sigma _{t}(b)) & = & f(t),\\
\phi (\sigma _{t}(b)a) & = & f(t+i),\quad t\in \reals .
\end{eqnarray*}
Fix \( P,Q\in B[X] \) and let \( f \) be as above, so that
\begin{eqnarray*}
\phi (\sigma _{t}^{\phi }(P)\xi Q) & = & f(t+i),\\
\phi (\xi Q\sigma ^{\phi }_{t}(P)) & = & f(t).
\end{eqnarray*}
Then
\begin{eqnarray*}
f(t) & = & \langle \xi ,Q\sigma ^{\phi }_{t}(P)\rangle \\
 & = & \langle Y,\partial _{X}(Q\sigma ^{\phi }_{t}(P))\rangle \\
 & = & \langle Y,\partial _{X}(Q)\sigma ^{\phi }_{t}(P)\rangle +\langle Y,Q\sigma ^{\hat{\phi }}_{t}(\partial _{X}(P))\rangle ,
\end{eqnarray*}
where in the last step we used the fact that \( \partial _{X} \) intertwines
\( \sigma _{t}^{\phi } \) and \( \sigma _{t}^{\hat{\phi }} \). Using the KMS-condition
for \( \hat{\phi } \), we then get
\begin{eqnarray*}
f(t+i) & = & \hat{\phi }(\sigma _{t}^{\hat{\phi }}(P)Y\partial _{X}(Q))+\hat{\phi }(\sigma _{t}^{\hat{\phi }}(\partial _{X}(P))YQ)\\
 & = & \hat{\phi }(\sigma _{t}^{\phi }(P)Y\partial _{X}(Q))+\hat{\phi }(\partial _{X}(\sigma _{t}^{\phi }(P))YQ).
\end{eqnarray*}
Since \( f(t+i)=\phi (\sigma _{t}^{\phi }(P)\xi Q), \) we get, setting \( t=0 \),
that
\[
\phi (P\xi Q)=\hat{\phi }(PY\partial _{X}(Q))+\hat{\phi }(\partial _{X}(P)YQ),\]
as claimed.
\end{proof}

\subsection{Conjugate variables as free Brownian gradients.}

As pointed out above, \( X+\sqrt{\varepsilon }Y \) is a natural free Brownian
motion, which is covariant with respect to the appropriate modular groups. The
following proposition shows that \( J_{\phi }(X:B) \) plays the role of the
free Brownian gradient of \( X \).

\begin{prop}
Assume that \( \xi =J_{\phi }(X:B) \) exists and belongs to \( M\subset L^{2}(M,\phi ) \).
Let \( P(Z_{1},\dots ,Z_{n}) \), be any non-commutative polynomial in \( n \)
variables \( Z_{1},\dots ,Z_{n} \), with coefficients from \( B \). Write
\( X_{t}=\sigma ^{\phi }_{t}(X) \), \( Y_{t}=\sigma ^{\hat{\phi }}_{t}(Y) \),
\( \xi _{t}=\sigma _{t}^{\phi }(\xi ) \), \( t\in \reals  \). 

Then for all \( t_{1},\dots ,t_{n}\in \reals  \), we have
\begin{eqnarray*}
\hat{\phi }(P(X_{t_{1}}+\sqrt{\varepsilon }Y_{t_{1}},\ldots ,X_{t_{n}}+\sqrt{\varepsilon }Y_{t_{n}})) & = & \frac{1}{2}\phi (P(X_{t_{1}}+\varepsilon \xi _{t_{1}},\ldots ,X_{t_{n}}+\varepsilon \xi _{t_{n}})\\
 &  & +O(\varepsilon ^{2}).
\end{eqnarray*}

\end{prop}
\begin{proof}
We may assume, by linearity, that \( P \) is a monomial, i.e., \( P(Z_{1},\dots ,Z_{n})=b_{0}Z_{1}b_{1}\cdots b_{n-1}Z_{n}b_{n} \),
for \( b_{j}\in B \). In this case, we have
\begin{eqnarray*}
\hat{\phi }(b_{0}(X_{t_{1}}+\sqrt{\varepsilon }Y_{t_{1}})b_{1}\cdots b_{n}) & = & \phi (P(X_{t_{1}},\ldots ,X_{t_{n}}))+O(\varepsilon ^{2})\\
 &  & +\varepsilon \sum _{k<l}\hat{\phi }(b_{0}X_{t_{1}}\cdots b_{l}Y_{t_{l+1}}b_{l+1}X_{t_{l+2}}\cdots b_{k}Y_{t_{k+1}}\cdot \\
 &  & \quad \cdot b_{k+1}X_{t_{k+2}}\cdots X_{t_{n}}b_{n})\\
 & = & \phi (P(X_{t_{1}},\ldots ,X_{t_{n}}))+O(\varepsilon ^{2})\\
 &  & +\frac{1}{2}\varepsilon \sum _{l}\hat{\phi }(b_{0}X_{t_{1}}\cdots X_{t_{l}}b_{l}Y_{t_{l+1}}\partial _{X_{t_{l+1}}}(b_{l+1}X_{t_{l+2}}\cdots X_{t_{n}}b_{n}))\\
 &  & +\frac{1}{2}\varepsilon \sum _{k}\hat{\phi }(\partial _{X_{t_{k+1}}}(b_{0}X_{t_{1}}\cdots X_{t_{k}}b_{k})Y_{t_{k}}b_{k+1}X_{t_{k+2}}\cdots X_{t_{n}}b_{n})\\
 & = & \phi (P(X_{t_{1}},\ldots ,X_{t_{n}}))+O(\varepsilon ^{2})\\
 &  & +\frac{1}{2}\sum _{k}\phi (b_{0}X_{t_{1}}\cdots X_{t_{k}}b_{k}\varepsilon \xi _{t_{k}}b_{k+1}X_{t_{k+2}}\cdots X_{t_{n}}b_{n}),
\end{eqnarray*}
the last equality by Lemma \ref{lemma:xileftright}. This implies the statement
of the Lemma.
\end{proof}

\subsection{Examples of conjugate variables.}

\subsubsection{Tracial case.}

We have seen before that if \( \phi  \) is a trace, then the definition of
\( J_{\phi }(X:B) \) coincides with the definition of conjugate variables
given by Voiculescu, up to a constant (which has to do with the fact that we
choose \( Y \) so that \( \Vert Y\Vert _{L^{2}(\hat{\phi })}=\Vert X\Vert _{L^{2}(\phi )} \),
and not \( 1 \)). In particular,
\[
J_{\phi }(X:B)=J(X:B)\cdot \frac{1}{\Vert X\Vert _{L^{2}(\phi )}^{2}},\quad \textrm{if }\phi \, \textrm{is a trace}.\]

\subsubsection{Free quasi-free states.}

Let \( \mu  \) be a positive finite Borel measure on \( \reals  \), and let
\( \mathcal{H}_{\reals } \) be the real Hilbert space \( L^{2}(\reals ,\reals ,\mu ) \)
of \( \mu  \) -square-integrable real-valued functions. Denote by \( U_{t} \)
the representation of \( \reals  \) on \( \mathcal{H}_{\reals } \), given
by 
\[
(U_{t}f)(x)=e^{2\pi itx}f(x),\quad x,t\in \reals .\]
Let \( h \) denote the vector \( 1\in \mathcal{H}_{\reals } \), and consider
\[
M=\Gamma (\mathcal{H}_{\reals },U_{t})'',\quad \phi =\phi _{U},\quad X=s(h)\in M\]
 (see \cite{shlyakht:quasifree:big} for definitions and notation). 

Then \( X=J_{\phi }(X:\mathbb {C}) \). Indeed, set \( X_{t}=\sigma ^{\phi }_{t}(X)=s(U_{t}h) \);
then we have
\begin{eqnarray*}
\phi (X\cdot X_{t_{1}}\cdots X_{t_{n}}) & = & \sum _{k}\phi (XX_{t_{k}})\phi (X_{t_{1}}\cdots X_{t_{k-1}})\cdot \phi (X_{t_{k+1}}\cdots X_{t_{n}})\\
 & = & \sum _{k}\hat{\phi }(Y\sigma ^{\hat{\phi }}_{t_{k}}(Y))\phi (X_{t_{1}}\cdots X_{t_{k-1}})\cdot \phi (X_{t_{k+1}}\cdots X_{t_{n}})\\
 & = & \sum _{k}\hat{\phi }(YX_{t_{1}}\cdots X_{t_{k-1}}YX_{t_{k+1}}\cdots X_{t_{n}})\\
 & = & \hat{\phi }(Y\partial _{X}(X_{t_{1}}\cdots X_{t_{n}})),
\end{eqnarray*}
so that \( X \) satisfies the defining property of \( J_{\phi }(X:\mathbb {C}) \),
and hence \( J_{\phi }(X:\mathbb {C}) \) exists and equals \( X \).

\subsection{Free Fisher information.}

Following \cite{dvv:entropy5}, we define the free Fisher information \( \Phi _{\phi }^{*}(X:B) \)
to be
\[
\Phi _{\phi }^{*}(X:B)=\Vert J_{\phi }(X:B)\Vert _{2}^{2}\cdot \Vert X\Vert _{2}^{-2}\]
(the extra factor \( \Vert X\Vert _{2}^{-2} \) comes from the fact that \( \partial _{X}(X) \)
does not have unit norm in our definition). For several variables, we set
\[
\Phi _{\phi }^{*}(X_{1},\dots ,X_{n})=\sum \Phi _{\phi }^{*}(X_{i}:W^{*}(\sigma _{t_{1}}^{\phi }(X_{1}),\ldots ,\hat{X}_{i},\ldots ,\sigma _{t_{n}}^{\phi }(X_{n}):t_{1},\ldots ,t_{n}\in \mathbb {R}))\]
(here \( \hat{X}_{i} \) means that \( X_{i} \) is omitted).

\section{Free Fisher Information Relative to the Core.}

Recall \cite{stratila:modular} that if \( (M,\phi ) \) is as above, its core
is defined to be the von Neumann algebra crossed product \( P=M\rtimes _{\sigma ^{\phi }}\reals  \).
There is a canonical inclusion \( M\subset P \), and \( P \) is densely spanned
by elements of the form
\[
mU_{t},\quad t\in \reals ,\]
where \( m\in M \), and \( U_{t} \) satisfy \( U_{t}mU_{t}^{*}=\sigma _{t}^{\phi }(m) \).
The elements \( U_{t}:t\in \reals  \) generate a copy of the group von Neumann
algebra \( L(\reals )\subset P \); the map
\[
E^{\phi }:mU_{t}\mapsto \phi (m)U_{t},\quad m\in M,t\in \reals \]
extends to a normal conditional expectation from \( P \) onto \( L(\reals ) \).

For \( X\in M \) self-adjoint, define the completely positive map \( \eta _{X}:L(\reals )\to L(\reals ) \)
by
\[
\eta _{X}(g)=E^{\phi }(XgX),\quad g\in L(\reals ).\]
Identify \( L(\reals ) \) with \( L^{\infty }(\reals ) \) via Fourier transform.
For each \( t\in \reals  \), set
\[
\eta (t)=\langle X,\sigma _{t}^{\phi }(X)\rangle =E^{\phi }(XU_{t}X).\]
Then \( \eta _{X}(f)=\hat{\eta }*f \), if \( f\in L^{\infty }(\reals )\cong L(\reals ) \);
here \( \hat{\eta } \) denotes Fourier transform. 

Define on \( P \) an \( L(\reals ) \)-valued inner product
\[
\langle a,b\rangle _{L(\reals )}=E^{\phi }(a^{*}b),\quad a,b\in P.\]
 Denote by \( L^{2}(P,E^{\phi }) \) the \( L(\reals ) \)-Hilbert bimodule
arising from the completion of \( P \) with respect to the norm induced by
this inner product. Note that the restriction of \( \langle \cdot ,\cdot \rangle _{L(\reals )} \)
to \( M\subset P \) is valued in the complex field, and coincides with the
inner product \( \langle a,b\rangle =\phi (a^{*}b) \) on \( L^{2}(M) \).

Denote by \( \langle \cdot ,\cdot \rangle _{\eta } \) the \( L(\reals ) \)-valued
inner product on \( P\otimes P \) (algebraic tensor product) given by
\[
\langle a\otimes b,a'\otimes b'\rangle _{\eta }=E^{\phi }(b^{*}\eta (E^{\phi }(a^{*}a')b')),\quad a,a',b,b'\in P.\]
 Denote by \( 1\otimes _{\eta }1 \) the vector \( 1\otimes 1\in P\otimes P \). 

Let \( \delta _{X}:B[X]\cdot L(\reals )\to P\otimes P \) be given by
\[
\delta _{X}(X)=1\otimes _{\eta }1,\quad \delta _{X}(B\cdot L(\reals ))=0\]
and the fact that \( \delta _{X} \) is a derivation.

\begin{thm}
Let \( (M,\phi ) \) be as above, and let \( P \) be its core. Let \( i:L^{2}(M,\phi )\to L^{2}(P,E^{\phi }) \)
be the extension of the inclusion of \( M\subset P \). Then \( \zeta =i(J_{\phi }(X:B)) \)
satisfies
\begin{equation}
\label{eqn:defofzeta}
\langle \zeta ,Q\rangle _{L(\reals )}=\langle 1\otimes _{\eta }1,\delta _{X}(Q)\rangle _{\eta _{X}}
\end{equation}
for all \( Q\in B[X]\vee L(\reals ) \). Conversely, if there exists a vector
\( \zeta \in L^{2}(P,E^{\phi }) \), so that (\ref{eqn:defofzeta}) is satisfied,
then \( J_{\phi }(X:B) \) exists and \( \zeta =i(J_{\phi }(X:B)) \).
\end{thm}
\begin{proof}
Assume first that \( J_{\phi }(X:B) \) exists. Set \( \zeta =i(J_{\phi }(X:B)) \).
We must verify that (\ref{eqn:defofzeta}) holds. By linearity, and the fact
that \( L(\reals )BL(\reals )\subset BL(\reals ) \), it is sufficient to consider
the case when \( Q=b_{0}U^{s_{1}}X_{t_{1}}b_{1}U^{s_{2}}\cdots X_{t_{n}}b_{n}U^{s_{n}} \),
with \( b_{j}\in B \) and \( X_{t}=\sigma _{t}^{\phi }(X) \). Then \( Q=P\cdot U^{r} \),
where \( r=\sum s_{j} \) , and \( P=b_{0}X_{t_{1}'}b_{1}'\cdots b_{n}' \),
with \( b_{j}'=\sigma ^{\phi }_{s_{j-1}}\circ \cdots \sigma _{s_{1}}^{\phi }(b_{j}) \),
\( t_{j}'=s_{j-1}+\cdots +s_{1}+t_{j} \). Note that for \( x,y,x',y'\in P \),
\( \langle x\otimes y,x'\otimes Uy'\rangle _{\eta }=\langle x\otimes y,x'U\otimes y\rangle _{\eta } \),
and \( \langle x\otimes y,U^{r}(x'\otimes y')U^{s}\rangle =\langle U^{-r}x\otimes yU^{-s}\rangle  \).
Using this , we get

\begin{eqnarray*}
\langle \zeta ,Q\rangle _{L(\reals )} & = & \langle \zeta ,P\rangle _{L(\reals )}g\\
 & = & \langle \zeta ,P\rangle _{L^{2}(M,\phi )}\cdot U^{r}\\
 & = & \hat{\phi }(Y\partial _{X}(P))\cdot U^{r}\\
 & = & \sum _{j}\phi (b'_{0}X_{t'_{1}}\cdots X_{t'_{j}}b'_{j})\phi (b'_{j+1}X_{t'_{j+2}}\cdots X_{t'_{n}}b'_{n})\cdot \hat{\phi }(YY_{t_{j}})U^{r}\\
 & = & \sum _{j}\phi (b_{0}'X_{t'_{1}}\cdots X_{t'_{j}}b_{j}')\phi (b'_{j+1}X_{t'_{j+2}}\cdots X_{t'_{n}}b'_{n})\cdot \phi (XU^{t'_{j}}XU^{-t'_{j}})U^{r}\\
 & = & \sum _{j}\phi (b'_{0}X_{t'_{1}}\cdots X_{t'_{j}}b'_{j})\phi (b'_{j+1}X_{t'_{j+2}}\cdots X_{t'_{n}}b'_{n})\cdot E^{\phi }(XU^{t'_{j}}XU^{-t'_{j}}U^{r})\\
 & = & \sum _{j}\phi (b'_{0}X_{t'_{1}}\cdots X_{t'_{j}}b'_{j})\phi (b'_{j+1}X_{t'_{j+2}}\cdots X_{t'_{n}}b'_{n})\cdot \eta _{X}(U^{t'_{j}})U^{r-t'_{j}}\\
 & = & \sum _{j}\eta _{X}\circ E^{\phi }(b'_{0}X_{t'_{1}}\cdots X_{t'_{j}}b'_{j}U^{t'_{j}})\cdot E^{\phi }(U^{r-t'_{j}}b'_{j+1}X_{t'_{j+2}}\cdots X_{t'_{n}}b'_{n})\\
 & = & \sum _{j}\langle 1\otimes _{\eta }1,b'_{0}X_{t'_{1}}\cdots X_{t'_{j}}b'_{j}U^{t'_{j}}\otimes U^{r-t_{j}}b'_{j+1}X_{t'_{j+2}}\cdots X_{t'_{n}}b'_{n}\rangle _{\eta }\\
 & = & \sum _{j}\langle 1\otimes _{\eta }1,b_{0}U^{s_{1}}X_{t_{1}}\cdots X_{t_{j}}b_{j}U^{s_{j}}U^{t'_{j}}\otimes U^{r-t_{j}}U^{s_{j+1}}b_{j+1}X_{t_{j+2}}\cdots X_{t_{n}}b_{n}U^{s_{n}}U^{-r}\rangle _{\eta }\\
 & = & \sum _{j}\langle 1\otimes _{\eta }1,b_{0}U^{s_{1}}X_{t_{1}}\cdots X_{t_{j}}b_{j}U^{s_{j}}\otimes U^{s_{j+1}}b_{j+1}X_{t_{j+2}}\cdots X_{t_{n}}b_{n}U^{s_{n}}\rangle _{\eta }\\
 & = & \langle 1\otimes _{\eta }1,\delta _{X}(Q)\rangle _{\eta }.
\end{eqnarray*}

Conversely, assume that \( \zeta  \) satisfying (\ref{eqn:defofzeta}) exists.
Since the argument above is reversible, it is sufficient to prove that \( \zeta  \)
is in the image of \( i:L^{2}(M)\to L^{2}(P,E^{\phi }) \). Let \( \theta _{t} \)
be the dual action of \( \reals  \) on \( P \), given by \( \theta _{t}(U_{s})=\exp (2\pi ist) \),
\( \theta _{t}(m)=m \), \( m\in M \). It is sufficient to prove that \( \theta _{t}(\zeta )=\zeta  \),
since \( i(L^{2}(M)) \) consists precisely of those vectors, which are left
fixed by \( \theta  \). It is sufficient to prove that \( \theta _{s}(E^{\phi }(\zeta mU_{t}))=\exp (2\pi ist) \)
if \( m\in M \). Since \( \zeta  \) is assumed to be in the closure of \( B[X]\vee L(\reals ) \),
it is sufficient to check this for \( m\in B[X] \). But then by (\ref{eqn:defofzeta}),
\begin{eqnarray*}
E^{\phi }(\zeta mU_{t}) & = & \langle \zeta ,mU_{t}\rangle _{L(\reals )}\\
 & = & \langle 1\otimes _{\eta }1,\delta _{X}(mU)\rangle _{\eta }\\
 & = & \langle 1\otimes _{\eta }1,\delta _{X}(m)\rangle _{\eta }U\\
 & \in  & M\cdot U_{t},
\end{eqnarray*}
which gives the desired result, since \( \theta _{s} \) acts trivially on \( M \).
\end{proof}
Note that (\ref{eqn:defofzeta}) means that \( \zeta  \) is equal to \( J(X:B\vee \mathbb {R},\eta ) \)
in the notation of \cite{shlyakht:cpentropy}. (This is strictly speaking incorrect,
since the setting of \cite{shlyakht:cpentropy} presumes the existence of a
finite trace on \( B[X]\vee L(\mathbb {R}) \); however, it is not hard to check
that the arguments in \cite{shlyakht:cpentropy} go through also in the case
of a semifinite trace, which exists in our case).

This fact has many consequences for the conjugate variables \( J_{\phi }(X:B) \),
coming from the properties of \( J(X:B\vee L(\mathbb {R}),\eta ) \). Note in
particular that if \( X \) is free from \( B \) with amalgamation over \( D\subset B \)
with respect to some conditional expectation \( E:B\to D \), and \( E \) is
\( \phi  \)-preserving, then \( X \) is free from \( B\vee L(\mathbb {R}) \)
with amalgamation over \( D\vee L(\mathbb {R}) \) (see \cite{ueda:SUqActionFullFactor},
\cite{shlyakht:semicirc}). We record this as

\begin{thm}
Assume that \( E:B\to D \) is a \( \phi  \)-preserving conditional expectation.
If \( X \) is free from \( B \) over \( D \), then
\[
J_{\phi }(X:B)=J_{\phi }(X:D).\]

\end{thm}
In a similar way, one can generalize to \( J_{\phi }(X:B) \) all the properties
of the conjugate variable \( J(X:B,\eta ) \) proved in \cite{shlyakht:cpentropy}. 

Reformulating gives the following properties of \( \Phi _{\phi } \), which
we list for reader's convenience, since they are needed in the rest of the paper:

\begin{thm}
Let \( \phi  \) be a normal faithful state on \( M \), \( B\subset M \) be
globally fixed by the modular group (i.e., \( \sigma _{t}^{\phi }(B)=B \) for
all \( t \)), and \( X_{i}\in M \). Then:\\
(a) \( \Phi ^{*}_{\phi }(\lambda X_{1},\ldots ,\lambda X_{n}:B)=\lambda ^{-2}\Phi ^{*}_{\phi }(X_{1},\ldots ,X_{n}:B) \)
for all \( \lambda \in \mathbb {R}\setminus \{0\} \)\\
(b) If \( B\subset A\subset M \) and \( A \) is globally fixed by \( \sigma ^{\phi } \),
then \( \Phi ^{*}_{\phi }(X_{1},\ldots ,X_{n}:A)\geq \Phi ^{*}(X_{1},\ldots ,X_{n}:B) \).\\
(c) If \( C\subset M \) is globally fixed by \( \sigma ^{\phi } \), and \( W^{*}(X_{1},\ldots ,X_{n}) \)
and \( B \) are free with amalgamation over \( C \) (with respect to the unique
\( \phi  \)-preserving conditional expectation from \( M \) onto \( C \)),
then \( \Phi ^{*}_{\phi }(X_{1},\dots ,X_{n}:B\vee C)=\Phi ^{*}_{\phi }(X_{1},\dots ,X_{n}:B) \).\\
(d) If \( Y_{i}\in M \) are self-adjoint, \( D\subset B \), \( D\subset C \)
subalgebras of \( M \), which are globally fixed by \( \sigma ^{\phi } \),
and \( B[X_{1},\ldots ,X_{n}] \) is free from \( C[X_{1},\ldots ,X_{n}] \)
over \( D \) (with respect to the unique \( \phi  \)-preserving conditional
expectation from \( M \) onto \( D \)), then \( \Phi ^{*}_{\phi }(X_{1},\ldots ,X_{n},Y_{1},\ldots ,Y_{m}:B\vee C)=\Phi ^{*}_{\phi }(X_{1},\ldots ,X_{n}:B)+\Phi ^{*}_{\phi }(Y_{1},\ldots ,Y_{n}:C) \).\\
(e) \( \Phi ^{*}_{\phi }(X_{1},\ldots ,X_{n},Y_{1},\ldots ,Y_{n}:B)\geq \Phi _{\phi }^{*}(X_{1},\ldots ,X_{n}:B)+\Phi ^{*}_{\phi }(Y_{1},\ldots ,Y_{m}:B) \).
\\
(f) \( \Phi _{\phi }^{*}(X_{1},\ldots ,X_{n}:B)\cdot \phi (\sum X_{i}^{*}X_{i})^{2}\geq n^{2} \).
Equality holds iff \( \{\sigma _{t_{1}}^{\phi }(X_{1}),\dots ,\sigma _{t_{n}}^{\phi }(X_{n}):t_{1},\dots ,t_{n}\in \mathbb {R}\} \)
have the same distribution as the semicircular family \( \{\kappa s(\sigma _{t_{1}}^{\phi }(X_{1})),\dots ,\kappa s(\sigma _{t_{n}}^{\phi }(X_{n})):t_{1},\dots ,t_{n}\in \mathbb {R}\} \)
with respect to the free quasi-free state, \( \kappa >0 \).
\end{thm}
We mention that all of the statements in sections 3 and 4 of \cite{shlyakht:cpentropy}
remain valid for \( \Phi ^{*}_{\phi } \); we leave details to the reader.

One can also define and study free entropy \( \chi ^{*}_{\phi }(X_{1},\dots ,X_{n}) \)
by setting \( X_{i}^{\epsilon }=X_{i}+Y_{i} \) to be the free Brownian motion
described in the beginning of the paper, and letting
\[
\chi ^{*}_{\phi }(X_{1},\dots ,X_{n})=\frac{1}{2}\int _{0}^{\infty }\left( \frac{n}{1+t}-\Phi ^{*}_{\phi }(X_{1}^{\epsilon },\ldots ,X_{n}^{\epsilon })\right) dt.\]
The properties of \( \chi ^{*}(\cdots ,\eta ) \) once again generalize to \( \chi ^{*}_{\phi } \)
(compare section 8 of \cite{shlyakht:cpentropy}).

\section{States on a II\protect\( _{1}\protect \) factor.}

\subsection{\protect\( \Phi _{\phi }^{*}\protect \) vs. \protect\( \Phi _{\tau }^{*}\protect \).}

The following theorem is somewhat surprising, since it shows that \( \Phi _{\phi }^{*} \)
is identically infinite for most states \( \phi  \) on a II\( _{1} \) factor
(the analogy with classical Fisher information would instead suggest that \( \phi \mapsto \Phi ^{*}_{\phi } \)
would have some nice convexity properties). This, on the other hand, goes well
with the ``degenerate convexity'' property of the microstates free entropy
\( \chi  \) \cite{dvv:entropy3} (which is reflected in that it is identically
\( -\infty  \) on generators of any von Neumann algebra which more than one
unital trace).

\begin{thm}
\label{thrm:mustbetrace}Let \( M \) be a tracial von Neumann algebra, \( \phi  \)
a faithful normal state on \( M \), \( B\subset M \) a subalgebra so that
\( \sigma _{t}^{\phi }(B)=B \) for all \( t \), and \( X=X^{*}\in M \). Then
if \( J_{\phi }(X:B) \) exists, the modular group of \( \phi  \) must fix
\( X \).
\end{thm}
\begin{proof}
Let \( d\in M \) be a positive element, so that \( \phi (x)=\tau (dx) \),
where \( \tau  \) is a normal faithful trace on \( M \), and \( d \) is an
unbounded operator on \( L^{2}(M,\tau ) \), affiliated to \( M \). The modular
group of \( \phi  \) is then given by \( \sigma _{t}^{\phi }(x)=d^{it}xd^{-it} \),
\( x\in M \). Denoting by \( X_{t} \) the element \( \sigma _{t}^{\phi }(X) \),
we then get
\[
X=X_{0}=d^{-it}X_{t}d^{it},\quad t\in \mathbb {R}.\]
Consider
\[
\phi (X_{0}^{2})=\phi (J_{\phi }(X:B)\cdot X_{0})=\phi (J_{\phi }(X:B)d^{-it}X_{t}d^{it}).\]
Let \( b_{1} \) and \( b_{2} \) be two elements in the domain of \( \partial _{X} \),
so that \( b_{1}=b_{2}^{*} \). Then we get, writing \( Y_{t}=\sigma ^{\hat{\phi }}_{t}(Y) \):
\begin{eqnarray*}
\phi (J_{\phi }(X:B)b_{1}X_{t}b_{2}) & = & \hat{\phi }(Y_{0}b_{1}Y_{t}b_{2})\\
 &  & +\hat{\phi }(Y_{0}\partial _{X}(b_{1})X_{t}b_{2})+\hat{\phi }(Y_{0}b_{1}X_{t}\partial _{X}(b_{2}))\\
 & = & \phi (b_{1})\phi (b_{2})\hat{\phi }(Y_{0}Y_{t})\\
 &  & +\hat{\phi }(Y_{0}\partial _{X}(b_{1})X_{t}b_{2})+\hat{\phi }([\partial _{X}(b_{2}^{*})X_{t}b_{1}^{*}Y_{0}]^{*})\\
 & = & \phi (b_{1})\phi (b_{1}^{*})\hat{\phi }(Y_{0}Y_{t})\\
 &  & +\hat{\phi }(Y_{0}\partial _{X}(b_{1})X_{t}b_{1}^{*})+\hat{\phi }([\partial _{X}(b_{1})X_{t}b_{1}^{*}Y_{0}]^{*}).
\end{eqnarray*}
Now, for all \( m,n\in M \), we have
\[
\hat{\phi }(Y_{0}mY_{0}n)=\phi (m)\phi (n)=\hat{\phi }(mY_{0}nY_{0}),\]
so that
\begin{eqnarray*}
\hat{\phi }(Y_{0}\partial _{X}(b_{1})X_{t}b_{1}^{*})+\hat{\phi }([\partial _{X}(b_{1})X_{t}b_{1}^{*}Y_{0}]^{*}) & = & \hat{\phi }(Y_{0}\partial _{X}(b_{1})X_{t}b_{1}^{*})+\hat{\phi }([Y_{0}\partial _{X}(b_{1})X_{t}b_{1}^{*}]^{*})\\
 & = & \hat{\phi }(Y_{0}\partial _{X}(b_{1})X_{t}b_{1}^{*})+\overline{\hat{\phi }(Y_{0}[\partial _{X}(b_{1})X_{t}b_{1}^{*}])}\\
 & \in  & \mathbb {R}.
\end{eqnarray*}
It follows that
\[
\Im \phi (J_{\phi }(X:B)b_{1}X_{t}b_{1}^{*})=\Im \hat{\phi }(Y_{0}b_{1}Y_{t}b_{1}^{*}).\]
Now fix \( t\in \mathbb {R} \) and choose \( a_{n} \) in the domain of \( \partial _{X} \),
\( \Vert a_{n}\Vert \leq 1 \), so that
\[
a_{n}\to d^{it},\quad a_{n}^{*}\to d^{-it}\qquad \textrm{strongly}.\]
One can choose \( a_{n} \), for example, to be elements of the algebra \( B[X] \).
Then
\begin{eqnarray*}
0=\Im \phi (X_{0}^{2}) & = & \Im \phi (J_{\phi }(X:B)\cdot X_{0})\\
 & = & \Im \phi (J_{\phi }(X:B)d^{-it}X_{t}d^{it})\\
 & = & \lim _{n\to \infty }\Im \phi (J_{\phi }(X:B)a_{n}X_{t}a_{n}^{*})\\
 & = & \lim _{n\to \infty }\Im \hat{\phi }(Y_{0}a_{n}Y_{t}a_{n}^{*})\\
 & = & \lim _{n\to \infty }\Im (\phi (a_{n})\phi (a_{n}^{*})\hat{\phi }(Y_{0}Y_{t})\\
 & = & \lim _{n\to \infty }\phi (a_{n})\phi (a_{n}^{*})\Im \hat{\phi }(Y_{0}Y_{t})\\
 & = & \phi (d^{it})\phi (d^{-it})\Im \hat{\phi }(Y_{0}Y_{t}).
\end{eqnarray*}
Since \( \hat{\phi }(Y_{0}Y_{t})=\phi (X_{0}X_{t}) \), for \( t \) sufficiently
close to zero (so that \( \phi (d^{it})\neq 0 \)), we get that
\[
\phi (XX_{t})\in \mathbb {R}.\]
Thus
\begin{eqnarray*}
0 & = & \tau (dXX_{t})-\tau (d(XX_{t})^{*})\\
 & = & \tau (dXd^{it}Xd^{-it}-d^{it}Xd^{-it}Xd)\\
 & = & \tau ((dX-Xd)d^{it}Xd^{-it})\\
 & = & \tau ([d,X]d^{it}Xd^{-it}).
\end{eqnarray*}
Differentiating this in \( t \), and noting that \( (d/dt)_{t=0}(d^{it}Xd^{-it})=i[d,X] \)
gives
\[
i\tau ([d,X]^{2})=0.\]
Since \( [d,X] \) is anti-self-adjoint, this implies that \( \tau (|[d,X]|^{2})=0 \),
so that \( [d,X]=0 \), because \( \tau  \) is faithful. This means that \( \sigma _{t}^{\phi }(X)=X \)
for all \( t \).
\end{proof}
\begin{cor}
Suppose that \( X_{1},\dots ,X_{n} \) are self-adjoint generators of a II\( _{1} \)
factor \( M \). Let \( \phi  \) be a normal faithful state on \( M \), and
denote by \( \tau  \) the unique faithful normal trace on \( M \). Then \( \Phi ^{*}_{\phi }(X_{1},\dots ,X_{n})<+\infty  \)
implies that:\\
(1) \( \Phi ^{*}_{\tau }(X_{1},\dots ,X_{n})<\infty  \) and\\
(2) \( \phi  \) is a multiple of the trace \( \tau  \) on \( M \).
\end{cor}
\begin{proof}
Clearly, the second statement implies the first. To get the second statement,
write \( \phi (\cdot )=\tau (d\cdot ) \) and apply the theorem to conclude
that \( [d,X_{i}]=0 \). Since \( X_{1},\dots ,X_{n} \) generate \( M \),
\( d \) must be in the center of \( M \), which must consist of multiples
of identity, since \( M \) is a factor. But then \( d \) is a scalar multiple
of identity, so that \( \phi  \) and \( \tau  \) are proportional.
\end{proof}

\subsection{Factoriality.}

Voiculescu showed \cite{dvv:entropy3} that for his microstates entropy \( \chi  \)
the following implication holds:
\[
\chi (X_{1},\dots ,X_{n})>-\infty \Rightarrow W^{*}(X_{1},\dots ,X_{n})\textrm{ is a factor}.\]
In fact, the conclusion is stronger: not only is the center of \( W^{*}(X_{1},\dots ,X_{n}) \)
is trivial, but so is its asymptotic center. Unfortunately, we don't know if
the same implication holds for the non-microstates free entropy \( \chi ^{*} \)
introduced by Voiculescu in \cite{dvv:entropy5}, or even under the stronger
assumption that \( \Phi ^{*}(X_{1},\dots ,X_{n}) \) is finite. We prove a weaker
version of the assertion above for \( \Phi ^{*}=\Phi ^{*}_{\tau } \). We first
need a technical lemma:

\begin{lem}
\label{lemma:almostcommutation}Let \( \phi  \) be a normal faithful state
on \( M \). Let \( X\in M \) be self-adjoint and \( B\subset M \) be a subalgebra,
so that \( \sigma _{t}^{\phi }(B)=B \) for all \( t \). Assume that \( p\in B \)
is a self-adjoint projection, \( \phi (p)=\alpha  \), and so that \( \sigma _{t}^{\phi }(p)=p \)
for all \( p \). Assume that \( \Vert [X,p]\Vert _{2}<\delta  \). Then
\[
\Phi _{\phi }^{*}(X:B)>4\frac{\alpha ^{2}(1-\alpha )^{2}}{\delta ^{2}}.\]

\end{lem}
\begin{proof}
Let \( (A,\tau ) \) be a copy of \( L(\mathbb {F}_{2}) \), free from \( B[X] \).
Since \( \Phi ^{*}_{\phi }(X:B)=\Phi ^{*}_{\phi }(X:B\vee A) \), and since
the centralizer of \( \{\sigma _{t}^{\phi }(B)\}_{t}\vee A \) is a factor \cite{barnett},
we can find a projection \( q\in B\vee A \), which is fixed by the modular
group, and so that \( \Vert [X,q]\Vert _{2}<\delta  \), and \( \tau (q)=\beta =m/n \)
is rational and close to \( \alpha  \). We may moreover find a family of matrix
units \( e_{ij}\in B\vee A \), \( 1\leq i,j\leq n \), fixed by the modular
group, and so that
\begin{eqnarray*}
e_{ij}^{*}=e_{ji}, &  & e_{ij}e_{kl}=\delta _{jk}e_{il}\\
\tau (e_{ii})=\frac{1}{n}, &  & q=\sum _{i=1}^{m}e_{ii}.
\end{eqnarray*}
Denote by \( C \) the algebra generated in \( B\vee A \) by \( \{e_{ij}\}_{1\leq i,j\leq n} \).
Note that \( C\cong M_{n\times n} \), the algebra of \( n\times n \) matrices.
The restriction of \( \phi *\tau  \) to \( C \) is the usual matrix trace.
Then
\[
\Phi ^{*}_{\phi }(X:B)=\Phi ^{*}_{\phi *\tau }(X:B\vee A)\geq \Phi ^{*}_{\phi *\tau }(X:C).\]
Write \( X_{ij}=e_{1i}Xe_{j1} \). Then the inequality \( \Vert [X,q]\Vert _{2}<\delta  \)
implies that
\begin{eqnarray*}
\delta  & > & \Vert xq-qx\Vert _{2}\\
 & = & \Vert qxq+(1-q)xq-qxq-qx(1-q)\Vert _{2}\\
 & = & \Vert (1-q)xq-qx(1-q)\Vert _{2}\\
 & = & \sqrt{2}\cdot \Vert qx(1-q)\Vert _{2},
\end{eqnarray*}
since \( (1-q)xq \) and \( qx(1-q) \) are orthogonal. Hence
\[
\Vert qX(1-q)\Vert _{2}<\delta /\sqrt{2.}\]
It follows that
\[
\sum _{1\leq i\leq m,\, m<j\leq n}\phi (X_{ij}^{*}X_{ij})+\sum _{m<i\leq n,\, 1\leq i\leq n}\phi (X_{ij}^{*}X_{ij})<\delta ^{2}.\]
 Denote by \( \phi ' \) the state \( n(\phi *\tau )(e_{11}\cdot e_{11}) \)
on \( e_{11}W^{*}(X,C)e_{11} \). Then
\begin{eqnarray*}
\Phi ^{*}_{\phi '}(\{X_{ij}\}) & \geq  & \sum _{i,j}\Phi ^{*}_{\phi '}(\{X_{ij}\})\\
 & > & 2m(n-m)\frac{1}{n(\delta ^{2}/2m(n-m))}\\
 & = & \frac{(2m(n-m))^{2}}{n\delta ^{2}}\\
 & = & n^{3}4\frac{\beta ^{2}(1-\beta )^{2}}{\delta ^{2}}.
\end{eqnarray*}
Arguing exactly as in \cite[Proposition 4.1]{nss:entropy}, we get that
\[
\Phi _{\phi *\tau }^{*}(X:C)=\frac{1}{n^{3}}\Phi _{\phi '}^{*}(\{X_{ij}\})>4\frac{\beta ^{2}(1-\beta )^{2}}{\delta ^{2}}.\]
Since \( \beta  \) was a rational number, arbitrarily close to \( \alpha  \),
we get the desired estimate for \( \Phi ^{*}_{\phi }(X:B) \).
\end{proof}
\begin{thm}
\label{thrm:factorialitytracial}Assume that \( M \) is a von Neumann algebra
with a faithful normal trace \( \tau  \), and \( X_{i} \) are a family of
self-adjoint elements in \( M \), \( \Vert X_{i}\Vert =1 \). Assume that \( B_{i} \)
form an increasing sequence of subalgebras of \( M \), so that \( M=\overline{\cup B_{i}}^{w} \).
Assume further that for some normal faithful state \( \phi  \) on \( M \)
and for each \( j \),
\[
\liminf _{i}\Phi ^{*}_{\phi }(X_{i}:B_{j})<+\infty .\]
Then \( M \) is a factor.
\end{thm}
\begin{proof}
In view of Theorem \ref{thrm:mustbetrace}, we may assume that \( \phi  \)
is a trace, \( \tau  \). Assume that \( M \) is not a factor. Then there exists
a central projection \( p\in M \) of some trace \( \alpha =\tau (p) \), \( \alpha (1-\alpha )\neq 0 \).
Moreover, \( [p,X_{i}]=0 \) for all \( i \). Since \( B_{i} \) increase to
all of \( M \), given \( \delta >0 \), there is a large enough \( j \) and
a projection \( q\in B_{j} \), so that \( \Vert q-p\Vert _{2}<\delta /2 \).
Then for any \( k \),

\begin{eqnarray*}
\Vert [q,X_{k}]\Vert _{2} & = & \Vert qX_{k}-X_{k}q\Vert _{2}\\
 & = & \Vert (q-p)X_{k}-X_{k}(q-p)+pX_{k}-X_{k}p\Vert _{2}\\
 & \leq  & \Vert (q-p)X_{k}\Vert _{2}+\Vert X_{k}(q-p)\Vert _{2}+0\\
 & \leq  & 2\Vert (q-p)\Vert _{2}\Vert X_{k}\Vert \\
 & < & 2(\delta /2)=\delta .
\end{eqnarray*}
Now applying Lemma \ref{lemma:almostcommutation}, we deduce that for any \( i>j \),
\[
\Phi _{\tau }^{*}(X_{i}:B_{j})>4\frac{\alpha ^{2}(1-\alpha )^{2}}{\delta ^{2}}.\]
Hence \( \liminf _{i}\Phi _{\tau }^{*}(X_{i}:B_{j})>4\alpha ^{2}(1-\alpha )^{2}/\delta ^{2} \),
which is a contradiction, since \( \delta  \) was arbitrary.
\end{proof}
The hypothesis of the theorem is satisfied for some von Neumann algebras. For
example, let \( M=L(\mathbb {F}_{\infty }) \) generated by an infinite semicircular
family \( X_{i} \), \( i=1,2,3,\ldots  \). Then if \( B_{i}=W^{*}(X_{j}:j<i) \),
the assumptions of the theorem are satisfied. In fact if \( Y_{i} \) is \emph{any}
family of elements of a tracial von Neumann algebra, so that \( \Vert Y_{i}\Vert =1 \),
and \( X_{i} \) are a free semicircular family, then letting \( Z_{j}(\epsilon )=Y_{i}+\epsilon X_{i} \),
\( M_{\epsilon }=W^{*}(Z_{1}(\epsilon ),Z_{2}(\epsilon ),\ldots ) \) and \( B_{j}=W^{*}(Z_{1}(\epsilon ),\ldots ,Z_{j}(\epsilon )) \),
we see that \( M_{\epsilon } \) is a factor. In other words, generators of
an arbitrary tracial von Neumann algebra can be perturbed (in a certain representation
of this algebra) by an arbitrarily small amount \( \epsilon  \) in uniform
norm, to produce a II\( _{1} \) factor. Another way of putting it is to note
that the free Brownian motion \( \epsilon \mapsto Z_{j}(\epsilon ) \) started
at \( \{Y_{1},Y_{2},\ldots \} \) generates a factor at any time \( \epsilon >0 \).

\subsection{Factoriality in the non-tracial case. }

In a similar way, we get the following:

\begin{thm}
Assume that \( M \) is a von Neumann algebra with a faithful normal state \( \phi  \),
and \( X_{i} \) are a family of self-adjoint elements in \( M \), \( \Vert X_{i}\Vert =1 \).
Assume that \( B_{i} \) form an increasing sequence of subalgebras of \( M \),
\( \sigma _{t}^{\phi }(B_{i})=B_{i} \) for all \( t \) and \( i \), and assume
that \( M^{\phi }=\overline{\cup (B_{i}\cap M^{\phi })}^{w} \). Let \( R_{i} \)
be the operator of right multiplication by \( X_{i} \) densely defined on \( L^{2}(M,\phi ) \).
Assume that \( \sup _{i}\Vert R_{i}\Vert =C<+\infty  \). Assume further that
for each \( j \),
\[
\liminf _{i}\Phi ^{*}_{\phi }(X_{i}:B_{j})<+\infty .\]
Then \( M \) is a factor.
\end{thm}
\begin{proof}
Assume that \( M \) is not a factor. Then there exists a central projection
\( p\in M \), \( \alpha =\phi (p) \), \( \alpha (1-\alpha )\neq 0 \). Moreover,
\( [p,X_{i}]=0 \) for all \( i \). Since automatically \( p\in M^{\phi } \)
and \( B_{i}\cap M^{\phi } \) increase to all of \( M^{\phi } \), given \( \delta >0 \),
there is a large enough \( j \) and a projection \( q\in B_{j}\cap M^{\phi } \),
so that \( \Vert q-p\Vert _{2}<\delta /(1+C) \). Then for any \( k \),
\begin{eqnarray*}
\Vert [q,X_{k}]\Vert _{2} & = & \Vert qX_{k}-X_{k}q\Vert _{2}\\
 & = & \Vert (q-p)X_{k}-X_{k}(q-p)+pX_{k}-X_{k}p\Vert _{2}\\
 & \leq  & \Vert (q-p)X_{k}\Vert _{2}+\Vert X_{k}(q-p)\Vert _{2}+0\\
 & \leq  & \Vert R_{k}\Vert \Vert (q-p)\Vert _{2}+\Vert (q-p)\Vert _{2}\Vert X_{k}\Vert \\
 & < & (1+C)(\delta /(1+C))=\delta .
\end{eqnarray*}

\end{proof}
Note that the assumption on the norms of \( R_{i} \) is satisfied if each \( X_{i} \)
is analytic for \( \sigma _{t}^{\phi } \) and satisfies \( \sup _{k}\Vert \sigma _{i}^{\phi }(X_{k})\Vert =C<+\infty  \). 

\bibliographystyle{amsplain}

\providecommand{\bysame}{\leavevmode\hbox to3em{\hrulefill}\thinspace}

\end{document}